# Mathematics A Missing Factor in Assessing the Antiquity of *Vedāṅga Jyotiṣa*

Jaidev Dasgupta

While reading ancient texts one has to be cognizant of the assumptions made about the past. One has to ask: Are these assumptions valid? Are we projecting the present views into the past? A case in point is the dating of *Vedāṅga Jyotiṣa*. This note reports that due to mathematics being overlooked in the various attempts of assessing the time of composition of the text, the latter was dated much before the period when its mathematics was actually feasible.

*Jyotiṣa* (astronomy), along with *Śikṣā* (phonetics), *Chandas* (prosody), *Vyākaran* (grammar), *Nirukta* (etymology), and *Kalpa* (rituals) forms six limbs or *aṅgas* of the Vedic knowledge. These are not the parts of the Vedas but auxiliary disciplines the learning of which was essential for proper recitation, understanding, and practice of the Vedic texts. While the kalpas were about the rituals and how they were to be performed, the purpose of jyotiṣa was to determine suitable times for performing the rituals.

This ancient astronomy is known as the *Vedāṅga Jyotiṣa* (VJ hereafter) or *Jyotiṣa Vedāṅga* that has come down to us in two recensions – *Ṛg-Jyotiṣa* and *Yajur-Jyotiṣa*. While the former has 36 verses, 30 of which are in common with the latter that has 13 additional verses. In total there are 49 verses in VJ (Dikshit, 1969). For being a part of the Vedic knowledge and considered as the oldest text of mathematical astronomy on the Indian subcontinent used for calendar-making, its time of composition is of importance for the history of astronomy.

Multiple efforts have been made to date VJ. Different scholars have dated it variously. While Kuppanna Sastry and Dikshit have placed it sometime around 1400 to 1200 BCE (Sastry, 1984;

Dikshit, 1969) and Bárhaspatyah around 1098-99 BCE (Bárhaspatyah, 1907), Ohashi considers it composed between 600 and 400 BCE (Ohashi, 2016) and Pingree around 400 BCE (Pingree, 1981). Kim Plofker, in her book *Mathematics in India* (p. 39, Plofker, 2009), writes that "…a conclusive date for the text still eludes us. Linguistically it seems to belong to the post-Vedic, pre-Classical Sanskrit corpus, which would probably put it instead somewhere around the fifth or fourth century BCE."

The early dates are computed based on the time taken in the movement due to precession: (1) either of the winter solstice from its location mentioned in VJ at the first point of *Dhaniṣṭhā nakṣatra* to another known nakṣatra or asterism at a later point in time, such as in Varāhamihira's time in 530 CE at *Uttarāṣāḍhā* ¼ (p. 13, Sastry); (2) or of the shift in the longitude of alpha-Delphini, the junction star of Dhaniṣṭhā, in 1887 CE from the time in VJ (p. 87, Dikshit). Multiplying the difference in the longitudes of the two positions with the rate of precession (1° in 72 years) gives the number of years elapsed between the events. Subtracting the calculated number of years from the known year (530 CE or 1887 CE) gives the time of conjunction in VJ. Scholars have taken this time as that of the composition of VJ, which is an assumption. The composition could have occurred much after the time of conjunction.

Ohashi contests these early dates. According to him, a small error in observation with the naked eye can introduce a difference of hundreds of years. Moreover, he writes, "…it is difficult to believe that the Jyotiṣa-vedāṅga is earlier than the later Vedic period if we compare the developed system of astronomy in the Jyotiṣa-vedāṅga and the developing knowledge of astronomy in the Vedic literature." The 400 BCE estimation made by Pingree, on the other hand, is based on the observation that the language of VJ is post-Vedic: "it imitates Piṅgala's Chandaḥsūtra in using the final or first syllables of the names of the *nakṣatras* as their designations…" (Pingree, p. 10)

Surprisingly, none of these efforts asked about the feasibility of mathematics in VJ in the times they proposed. It was assumed, perhaps inadvertently, that the required knowledge of mathematics for astronomical calculations was already available then. Modern readers due to their familiarity with present mathematical systems and operations rarely ever stop to question whether this mathematics was feasible in the ancient times of VJ. Anachronistically, they project the present knowledge into antiquity.

As has been discussed in detail elsewhere (Dasgupta, 2023), mathematics in the Vedic period was concrete, i.e., it used objects such as bricks or, perhaps, sticks or tally marks on sticks in the absence of a symbolic representation of numbers. Mathematics using symbols could not have been possible until after the advent of *Brāhmi* script in the 3$^{rd}$ century BCE. Hence the kind of calculations in VJ, especially those involving fractions which depend on symbolic manipulation, had to wait for a long time.

Taking an example calculation from VJ: in a five-year yuga there are 67 sidereal months. During each of these months the moon traverses through 27 nakṣatras or asterisms. Hence it traverses through 1809 (67x27) nakṣatras in five years. In the same period, the moon also goes through 124 fortnights (2 fortnights in each of 62 synodic months in five years). Therefore, the moon covers 1809/124 (= 14.58871), or 14 73/124, asterisms per fortnight. Obviously, there is a fractional calculation involved here. As the moon keeps moving from one fortnight to another, the same number with its fractional part gets added repeatedly. As discussed below, such calculations were not possible until around the turn of the Common Era.

Aside from the above example, there are other calculations in VJ that depend on fractional arithmetic. For instance, the rule of three involves knowledge of fractions (ratios), their equality and the calculation of an unknown quantity, x, from three known quantities b, c, and d (x/b = c/d; therefore, x = b.c/d). This calculation, now known as proportions, is an algebraic procedure.

Therefore, the question arises: In the absence of numerals (because of no script), the place-value system, and zero, how was the above calculation actually possible in the Vedic and the early part of the post-Vedic period? It is well established that Brāhmi script made its first appearance in the 3$^{rd}$ century BCE in King Ashoka's edicts and took several centuries after that to further evolve and mature (Datta and Singh, 2004). Also, the decimal place-value notation and a symbol for zero appeared in Indian mathematics between 5$^{th}$ to 7$^{th}$ century CE (Katz 2009; Joseph 2016).

Other ancient texts, namely, *Arthaśāstra* (300 BCE-300 CE) and *Śulvasūtras* (800-300 BCE) do not offer much help in this regard either. While both display mostly the knowledge of unit fractions, it is practice bound, not using symbolic notations for calculation. In Arthaśāstra, unit fractions appear in the description of taxes, tolls, fee, surcharges and service charges, suggesting their use in practical life (Rangarajan, 1992). But that does not necessarily indicate the ability to do fractional calculations. In practice, for instance, 1/6$^{th}$ tax on agricultural production could be paid simply by separating the sixth measure for tax after every five measures by weight (of any unit) without any need for calculation. Furthermore, the measures, weights and balances used for different purposes were so well defined in whole numbers that a daily need for fractional calculation is doubtful.

Likewise, the *Śulvasūtras* – the manuals for constructing sacrificial altars – do talk about unit fractions. But the way numbers were tied to bricks and stones in altar-making in the Vedic period, geometry was also practice-bound and concrete – using bamboo, ropes and pegs (Dasgupta, 2023). Though by altering the length of a rope by a part of it one could draw a circle or a square of desired area, it is not the same as doing fractional calculations using symbolic numerals. A Jain text,

*Sthānāṅga Sūtra*, from around the turn of the era also mentions fractions (*kalāsavarna*) but offers no examples to show how they were computed.

Given these historical facts, the mathematics in the extant version of VJ seems inconsistent with the age it is portrayed in. Its post-Vedic origin also has to be reasonably later. Datta and Singh (2004) suggest that writing and reduction of fractions to lowest terms appeared in the Indian subcontinent around 200 CE, which was several centuries after the appearance of Brāhmi script during which symbols gradually evolved and developed. In the same duration, presumably, mathematical handling of fractions also developed. Modern research on evolution, teaching, and learning of mathematics shows that dealing with fractions is not easy (Devlin, 2000) unless one is familiar with multiplication, division, factorization, and the procedure of using the least common denominator. And without symbolic numerals, representing and performing these mathematical operations on fractions are incomprehensible. No wonder it took a while for methods to develop before *Āryabhaṭa* first reported the rules for handling fractions in his book *Āryabhaṭīya* around the end of the $5^{th}$ century. It may be assumed that, perhaps, some of these rules were known a century or so before him.

Therefore, the extant version of VJ has to be from the early centuries of this era. This timing also agrees with Pingree's observation of the language of VJ imitating Pingala. Since Pingala is dated around 200 BC (Bag, 1966; Singh 1986), composition of VJ has to be after that. Also, it must be noted that Arthaśāstra does not mention VJ. Though the text describes a five-year yuga, the equinoxes, the solstices, six seasons, twelve months, the synodic and sidereal months, the two fortnights, the intercalary months and when they are added, the astronomical conjunction of the sun, the moon, and the Dhaniṣṭhā asterism at winter solstice with which the five-year cycle begins in VJ is missing in it. Very likely, VJ was, therefore, created later by adding fractional calculations to a simpler, pre-existing calendar described in the Arthaśāstra, when appropriate mathematics became available. For the starting reference point of the five-year yuga, the conjunction at solstice reported in VJ was chosen whose knowledge existed either in an oral tradition that traversed over the preceding centuries, or was in a lesser-known text. With the lack of knowledge of the precession of the solstices, the sky was assumed to be fixed, making the calendar repeatable every five years. Winter solstices at Dhaniṣṭhā in VJ was a fixed event. Even Varahamihira in the 6$^{th}$ century had no idea of precession. Though, as mentioned above, he did notice a shift in the location of winter solstice (Sastry, 1984).

Turning back to the antiquity of VJ, while scholars such as Sastry and Dikshit assumed the date of composition same as that of the conjunction, Ohashi and Pingree were in the right direction in assessing VJ's date based on the knowledge in the text and its language. Though mathematics, a major factor in the text, got overlooked. Mathematics in VJ is more advanced than that of the age in which it is usually projected. Feasibility of its mathematics is crucial in deciding the age of the text. The confusion about its time, which still persists, is because of the oversight of this fact; it is the result of assuming existence of the present mathematical knowledge in the ancient past.

**Contact Information**: jaidevd101@gmail.com